\documentclass[10pt]{amsart}
\usepackage{thm-restate,amsmath,amssymb,enumerate,graphicx}

\newtheoremstyle{plainsl}%
	{\topsep}
	{\topsep}
	{\slshape} 
	{}
	{\normalfont\bfseries}
	{.}
	{ }
	{}
\swapnumbers

\newtheorem{theorem}{Theorem}[section]
\newtheorem{lemma}[theorem]{Lemma}
\newtheorem{corollary}[theorem]{Corollary}

\newcommand\cref[1]{Corollary~\ref{cor:#1}}

\renewcommand\proof{\noindent\textsl{Proof. }}
\newcommand\sqr[2]{{\vbox{\hrule height.#2pt
    \hbox{\vrule width.#2pt height#1pt \kern#1pt
        \vrule width.#2pt}\hrule height.#2pt}}}
\renewcommand\qed{%
	\ifmmode\eqno\sqr53
	\else\nolinebreak\ \hfill\sqr53\medbreak\fi}

\begin{document}

\sloppy
\title{Cores of Vertex Transitive Graphs}
\author[Roberson]{David E. Roberson \\
                  \\
                  {D}\lowercase{epartment of} C\lowercase{ombinatorics and} O\lowercase{ptimization}\\
                  U\lowercase{niversity of} W\lowercase{aterloo} \\
                  200 U\lowercase{niversity} A\lowercase{venue} W\lowercase{est} \\
                  W\lowercase{aterloo}, ON N2L 3G1,
                  C\lowercase{anada} \\
                  \\
                  E\lowercase{mail}: \texttt{\lowercase{droberso@math.uwaterloo.ca}}}


\date{\today}
\maketitle

\vspace{-.7cm}

\begin{abstract}
A core of a graph $X$ is a vertex minimal subgraph to which $X$ admits a homomorphism. Hahn and Tardif have shown that for vertex transitive graphs, the size of the core must divide the size of the graph.  This motivates the following question: when can the vertex set of a vertex transitive graph be partitioned into sets each of which induce a copy of its core? We show that normal Cayley graphs and vertex transitive graphs with cores half their size always admit such partitions. We also show that the vertex sets of vertex transitive graphs with cores less than half their size do not, in general, have such partitions.
\end{abstract}

\vspace{.05cm}

\section*{Introduction}\label{sec:intro}

A homomorphism from a graph $X$ to a graph $Y$ is map $\varphi: V(X) \rightarrow V(Y)$ which preserves adjacency, i.e.~$\varphi(x)$ is adjacent to $\varphi(y)$ whenever $x$ is adjacent to $y$. A \emph{core} of a graph $X$ is a vertex-minimal subgraph to which $X$ admits a homomorphism. It is a straightforward exercise to show that all cores of a graph are isomorphic, and therefore we sometimes refer to \emph{the} core of a graph. In this paper we are concerned with vertex transitive graphs, i.e.~graphs whose automorphism groups act transitively on their vertex sets.

This paper is motivated by a comment once made by my algebraic graph theory instructor who said that he does not know of any interesting examples of vertex transitive graphs that are neither cores nor have complete graphs as their cores. By this he meant that all other examples of vertex transitive graphs seem to simply be constructed by taking several copies of some vertex transitive core and then adding some edges between these copies. This remark motivates the question of when a vertex transitive graph can be partitioned into copies of its core, to which we are able to give a partial answer. In particular, we will show that such a partition exists for a subclass of Cayley graphs known as normal Cayley graphs, as well as for any vertex transitive graph whose core is half its size. As a bonus, the techniques developed for dealing with the normal Cayley graph case can be used to give an alternative proof of Theorem~\ref{vtxtranscores} below, which was originally proven in \cite{Tardif}.

\section{Normal Cayley Graphs}\label{normal}

In order to partition a graph into copies of its core, we must at least have that the size of the core divides the size of the graph. Fortunately, this is the case for all vertex transitive graphs, and in fact we have the following stronger statement from \cite{Tardif}.

\begin{restatable}[Hahn \& Tardif]{theorem}{fibres}\label{vtxtranscores}
If $X$ is a vertex transitive graph and $\varphi$ is an endomorphism of $X$ whose image is a core $Y$ of $X$, then all of the fibres $\varphi^{-1}(y)$, $y \in V(Y)$, have the same size and thus $|V(Y)|$ divides $|V(X)|$.
\end{restatable}

Note that an endomorphism is simply a homomorphism from a graph to itself, and a fibre of a homomorphism is the inverse image of a vertex. Typically, the term `fibre' is used to specifically refer to preimages of vertices contained in the image of a homomorphism, i.e.~nonempty preimages. Below, we will frequently use the notion of a \emph{retraction}, which is an endomorphism that is identity on its image. The image of a retraction is known as a \emph{retract}, and it is not difficult to see that a core of a graph is always a retract.\\

For a group $G$ and inverse closed subset $C \subseteq G\setminus\{1\}$, the \emph{Cayley graph} $X(G,C)$ is the graph with vertex set $G$ such that two vertices, $g$ and $h$, are adjacent if $g^{-1}h \in C$. The set $C$ is often referred to as the \emph{connection set}. A Cayley graph is \emph{normal} if its connection set is closed under conjugation by any group element, i.e.~ $g^{-1}Cg = C$ for all $g \in G$. When studying vertex transitive graphs, Cayley graphs are a natural class to consider since every vertex transitive graph is a retract of some Cayley graph~\cite{Sabidussi}.

For a Cayley graph $X = X(G,C)$ and fixed element $a \in G$, the map $f_a(g) = ag$ for all $g \in G$ is an automorphism of $X$ known as a left translation. The right translations, i.e.~ $g \mapsto ga$, are not necessarily automorphisms of $X$ since $(ga)^{-1}(ha) = a^{-1}(g^{-1}h)a$ is not necessarily an element of $C$ whenever $g^{-1}h \in C$. However, if $X$ is normal, then $a^{-1}(g^{-1}h)a \in C$ if and only if $g^{-1}h \in C$, and thus the right translations are automorphisms for normal Cayley graphs. In fact, an equivalent definition of normal Cayley graph is that the right translations are automorphisms.

To prove our main result concerning normal Cayley graphs, we will need the following two lemmas which will also allow us to give a new proof of Theorem~\ref{vtxtranscores}. The first is a simple yet useful lemma which applies to all graphs, not just those which are vertex transitive.

\begin{lemma}\label{orbitals}
Let $\varphi$ be an endomorphism of $X$ such that $\varphi(x) = \varphi(y)$ for distinct $x,y \in V(X)$, and let $w$ and $z$ be two vertices which appear in some core of $X$. Then there is no automorphism of $X$ which maps the pair $\{w,z\}$ to the pair $\{x,y\}$.
\end{lemma}
\proof
Suppose $\sigma$ is such an automorphism of $X$, and $\rho$ is a retraction onto a core of $X$ containing $w$ and $z$. Then the endomorphism $\varphi \circ \sigma \circ \rho$ has at least one fewer vertex in its image than is in the core of $X$, a contradiction.\qed

Applying this lemma to Cayley graphs, we obtain the following:

\begin{lemma}\label{Cayley}
Let $X = X(G,C)$ be a Cayley graph. If $\varphi$ is an endomorphism of $X$ whose image is a core $Y$ of $X$, and $y \in V(Y)$, then the sets $V(Y)a^{-1}$, for $a \in \varphi^{-1}(y)$, are mutually disjoint.
\end{lemma}
\proof
Suppose not. Then there exists distinct $a,b \in \varphi^{-1}(y)$ and distinct $c,d \in V(Y)$ such that $ca^{-1} = db^{-1}$. Consider the map $\sigma : V(X) \rightarrow V(X)$ given by $\sigma(x) = (ac^{-1})x$. Note that this is an automorphism of $X$ since it is a left translation. However, $\sigma(c) = a$ and $\sigma(d) = ac^{-1}d = b$ by the above. But this contradicts Lemma~\ref{orbitals}, since $\varphi(a) = y = \varphi(b)$ and both $c$ and $d$ appear in $Y$.\qed

We are now ready to show that normal Cayley graphs can always be partitioned into copies of their cores.

\begin{theorem}\label{corepartition}
Let $X$ be a normal Cayley graph and $Y$ be a core of $X$. Then there exists a partition $\{V_1,\ldots, V_k\}$ of $V(X)$ such that each $V_i$ induces a copy of $Y$.
\end{theorem}
\proof
Let $\varphi$ be a retraction from $X$ to $Y$. Further, let $A = \varphi^{-1}(y)$ for some $y \in V(Y)$. By Theorem~\ref{vtxtranscores}, $A$ has size $|V(X)|/|V(Y)|$, and thus $|A||V(Y)| = |V(X)|$. Combining this with the fact that the sets $V(Y)a^{-1}$, for $a \in A$, are mutually disjoint by Lemma~\ref{Cayley}, we see that these sets must in fact partition $V(X)$. Furthermore, since $X$ is normal, right translations are automorphisms, and therefore each set $V(Y)a^{-1}$, for $a \in A$, induces a copy of $Y$.\qed

Some of the ideas in the above proof are inspired by the proofs of Theorem~6.1.1 and Corollaries~6.1.2 and 6.1.3 in \cite{interestinggraphs}.

\section{Half-Sized Cores}

Given the result of Theorem~\ref{vtxtranscores}, it is natural to ask what can be said about vertex transitive graphs whose cores are half their size, since this is in some sense the simplest nontrivial case. It turns out that any such vertex transitive graph can be partitioned into two copies of its core. The proof uses the fact that the core of a vertex transitive graph is vertex transitive, a proof of which can be found in \cite{AGT} or \cite{Tardif}. We will use $x \sim_X y$ to denote that vertices $x$ and $y$ are adjacent in graph $X$, though we will often omit the subscript when it is clear from context. We use $X \cong Y$ to denote that $X$ and $Y$ are isomorphic graphs.

\begin{theorem}\label{halfsize}
Suppose $X$ is a vertex transitive graph with a core $X_1$ such that $|V(X_1)| = \frac{1}{2}|V(X)|$. Furthermore, let $\varphi: X \rightarrow X_1$ be a retraction onto $X_1$, and let $X_2$ be the subgraph of $X$ induced by the vertices $V(X) \setminus V(X_1)$. If $Y$ is the bipartite graph consisting of the edges of $X$ which have exactly one end in each of $V(X_1)$ and $V(X_2)$, then we have the following:
\begin{enumerate}
\item $X_1 \cong X_2$
\item $\varphi|_{_{X_2}}$ is an isomorphism from $X_2$ to $X_1$.
\item $Y$ is regular and all of its edges are of the form $\{x,\varphi(y)\}$ where $x \sim_{_{X_2}} y$.
\end{enumerate}
\end{theorem}

\proof
Since $X$ is vertex transitive, its core $X_1$ is also vertex transitive and therefore they are both regular. Let $d$ and $d_1$ be the degree of vertices in $X$ and $X_1$ respectively. This means that the $V(X_1)$ side of $Y$ is regular with degree $d - d_1$.

Since the fibres of $\varphi$ all have the same size, namely two, the restriction of $\varphi$ to $X_2$ is a bijection between the vertices of $X_2$ and $X_1$ that preserves adjacency and therefore $X_2$ is isomorphic to a spanning subgraph of $X_1$. So the degree of any vertex in $X_2$ is at most $d_1$ and thus the degree in $Y$ of a vertex in $V(X_2)$ is at least $d - d_1$. But of course this means that the degree in $Y$ of every vertex in $V(X_2)$ is exactly $d-d_1$ since the sum of the degrees on one side of a bipartite graph is equal to the sum of the degrees on the other side. Therefore $X_2$ is regular with degree $d_1$ and thus must be isomorphic to $X_1$, and furthermore, the restriction of $\varphi$ to $X_2$ is an isomoprhism from $X_2$ to $X_1$.

Note that we have already shown that $Y$ is regular with degree $d-d_1$. Hence, all that is left to show is that the edges of $Y$ have the appropriate form. Consider a vertex $x \in V(X_2)$ which is adjacent in $X$ to a vertex $y' \in V(X_1)$. Since $\varphi$ is a retraction, $\varphi(x) \sim \varphi(y') = y'$. However, since the restriction of $\varphi$ to $X_2$ is an isomorphism, if $y \in V(X_2)$ is such that $\varphi(y) = y'$, then $x \sim y$ in $X_2$.
\qed

The partition given in the above proof actually satisfies another interesting property; it is equitable. A partition $\{V_1,V_2,\ldots,V_k\}$ of the vertex set of a graph $X$ is said to be \emph{equitable} if for all $i,j \in [k]$, the number of neighbors a vertex of $V_i$ has in $V_j$ depends only on $i$ and $j$. Given such a partition, its \emph{quotient matrix} is the matrix whose $ij$-entry is the number of neighbors a vertex of $V_i$ has in $V_j$. It turns out that the eigenvalues of the quotient matrix of an equitable partition of a graph $X$, are all also eigenvalues of (the adjacency matrix of) $X$. From this and the above proof we obtain the following corollary.

\begin{corollary}
Let $X$ be a vertex transitive graph with a core $X_1$ such that $|V(X_1)| = \frac{1}{2}|V(X)|$, and let $X_2$ be the induced subgraph of $X$ with vertex set $V(X) \setminus V(X_1)$. Then $\{V(X_1),V(X_2)\}$ is an equitable partition of $X$. Furthermore, if $d$ is the valency of $X$ and $d_1$ is the valency of $X_1$, then $2d_1 - d \ge 0$ is an eigenvalue of of $X$.
\end{corollary}

\proof
The fact that $\{V(X_1),V(X_2)\}$ is an equitable partition is obvious from the above proof. Since it is an equitable partition, the eigenvalues of its quotient matrix are eigenvalues of $X$. As can also be easily seen from above, the quotient matrix of this equitable partition is
\[\left(\begin{array}{cc}
d_1 & d - d_1 \\
d - d_1 & d_1
\end{array}\right).\]
The eigenvalues of this matrix are $d$ and $2d_1 - d$.
\qed

A consequence of this corollary is that if a vertex transitive graph on $2p$ vertices for $p$ a prime is not bipartite and has no nonnegative integer eigenvalue of the same parity as its valency, then it is its own core.

It would be nice to be able to more precisely describe what the subgraph $Y$ in the above proof looks like. For instance, is it necessarily vertex transitive? A more precise description of $Y$ could lead to a complete characterization of when a graph $X$ is vertex transitive with a core of half its size.

\subsection{Attempts at Generalization}

In light of the above result, the obvious next step would be to attempt to generalize the argument to vertex transitive graphs with cores of smaller relative size. The first obstacle one runs into with this is that it is not obvious how to choose the partition. In the above, each fibre of the retraction onto the core had size two, and therefore one can choose the image of the retraction as one part of the partition, and the remaining vertices for the other part. In the more general case, there are many more possibilities for choosing the partition.

The second obstacle one runs into while trying to generalize the above result is that it is simply not true. For a counterexample, one can consider the line graph of $K_{2n}$, denoted $L(K_{2n})$, for $n \ge 3$. Since $K_{2n}$ has maximum degree $2n-1$ and can be $(2n-1)$-edge-colored, the graph $L(K_{2n})$ contains a clique of size $2n-1$ and can be $(2n-1)$-colored. Thus there are homomorphisms between $L(K_{2n})$ and $K_{2n-1}$ in both directions. Two such graphs are said to be \emph{homomorphically equivalent}, which implies that they have isomorphic cores, in this case $K_{2n-1}$. However, since $n \ge 3$, it is not hard to see that any clique of size $2n-1$ in $L(K_{2n})$ corresponds to the set of edges incident to a particular vertex in $K_{2n}$, and thus any two such cliques must share a vertex. Therefore, $L(K_{2n})$ cannot be partitioned into copies of its core. Since $|V(K_{2n-1})| = \frac{1}{n}|V(L(K_{2n}))|$, these graphs provide counterexamples to generalizing the above to any relative size smaller than $\frac{1}{2}$.

\subsection{Arc Transitive Graphs}

We say that a graph $X$ is \emph{arc transitive} if for any two ordered pairs of vertices $(x,y)$ and $(x',y')$ such that $x \sim y$ and $x' \sim y'$, there exists an automorphism, $\sigma$, of $X$ such that $\sigma(x) = x'$ and $\sigma(y) = y'$. It is easy to see that any arc transitive graph (without isolated vertices) is vertex transitive as well. As with vertex transitive graphs, the core of an arc transitive graph must be arc transitive as well. Furthermore, by a result in \cite{AGT}, the valency of an arc transitive graph is always divisible by the valency of its core. If we assume that the graph $X$ in Theorem`\ref{halfsize} is also arc transitive, then we can obtain even more specific results about its structure.

First, note that the graph $Y$ in the proof of Theorem~\ref{halfsize} has degree at most $d_1$ and therefore $X$ has degree $d \le 2d_1$. Since $d_1$ must divide $d$, this implies that either $d = d_1$ or $d = 2d_1$. In the former case $X$ is simply two disjoint copies of $X_1$, and in the latter case $X \cong X_1\left[\overline{K_2}\right]$. Therefore, an arc transitive graph has a core half its size if and only if it is either two disjoint copies of an arc transitive core, or is the lexicographic product of an arc transitive core and $\overline{K_2}$.

\section{Alternate Proof of Theorem~\ref{vtxtranscores}}

Here we present a new proof of Theorem~\ref{vtxtranscores} which uses the techniques of section~\ref{normal}. The proof uses what is known as the lexicographic product of a graph $X$ with a graph $Y$. This product is denoted $X[Y]$, and is the graph with vertex set $V(X) \times V(Y)$ and $(x,y) \sim (x',y')$ if either ($x = x'$ and $y \sim y'$) or $x \sim x'$. The special case of $Y$ being empty, i.e.~$X\left[\overline{K_m}\right]$, is known as a \emph{multiple} of $X$. This ``multiple" of $X$ is obtained by replacing every vertex of $X$ with an independent set of size $m$ and having all possible edges between two such independent sets corresponding to adjacent vertices of $X$.

Clearly, any multiple of $X$ contains a subgraph isomorphic to $X$, and the map which takes each independent set corresponding to a vertex of $X$ to that vertex of $X$ is a homomorphism to $X$. Therefore, $X$ is homomorphically equivalent to any of its multiples and thus they all have the same core. 

A well-known theorem by Sabidussi~\cite{Sabidussi} states that if $X$ is vertex transitive, then some multiple of $X$ is a Cayley graph. It is this theorem that allows us to give an alternative proof to Theorem~\ref{vtxtranscores}.

\fibres*
\proof
We first prove it for Cayley graphs. Suppose that $X = X(G,C)$ is a Cayley graph, and $\varphi$ is an endomorphism onto a core $Y$ of $X$. Then for any $y \in V(Y)$, the sets $V(Y)a^{-1}$, for $a \in \varphi^{-1}(y)$, are mutually disjoint by Lemma~\ref{Cayley}. This implies that $|\varphi^{-1}(y)| \le |V(X)|/|V(Y)|$ is true for all $y \in V(Y)$. However, the average size of a fibre of $\varphi$ is clearly $|V(X)|/|V(Y)|$, and thus we must have equality in the above inequality for all $y \in V(Y)$.

Now suppose that $X$ is a vertex transitive graph and $\varphi$ is an endomorphism whose image is a core $Y$ of $X$. By the result of Sabidussi mentioned above, there exists $m \in \mathbb{N}$ such that $Z = X\left[\overline{K_m}\right]$ is a Cayley graph. The vertices of $Z$ are of the form $(x,i)$ for $x \in V(X)$ and $i \in [m]$. The map $\rho$ given by $\rho(x,i) = (x,1)$ is easily seen to be a retraction onto a subgraph $X'$ of $Z$ isomorphic to $X$. If we define a map $\hat{\varphi} : V(X') \rightarrow V(X')$ by $\hat{\varphi}(x,1) = (\varphi(x),1)$, then $\hat{\varphi} \circ \rho$ is an endomorphism of $Z$ onto a subgraph $Y'$ isomorphic to $Y$, and is thus an endomorphism onto a core of $Z$. Therefore, by the first part of the proof, the fibres of $\hat{\varphi} \circ \rho$ all have the same size. However, each fibre of $\hat{\varphi} \circ \rho$ clearly has size $m$ times the size of the corresponding fibre of $\varphi$, and therefore all fibres of $\varphi$ must have the same size.\qed

\section*{Concluding Remarks}

This paper is a step towards fully describing the structure of vertex transitive graphs in terms of their cores. In the special cases of normal Cayley graphs and vertex transitive graphs with cores half their size, our results hint torward a product structure for these classes of graphs. The examples of $L(K_{2n})$ show us that not all vertex transitive graphs fit this description, and thus some more general result would be required to describe all vertex transitive graphs in terms of their cores. Currently, we do not know of any vertex transitive graphs which neither can be partitioned into copies of their core, nor have a complete graph as their core, and so it is an interesting question as to whether such graphs exist.

Another question of interest is when a the \emph{edge} set of a vertex transitive graph can be partitioned into copies of its core. A necessary condition for this is that the valency of the core must divide the valency of the graph, and thus arc transitive graphs are a natural class of graphs to consider with respect to this question. Note that the counterexamples for vertex partitioning, the graphs $L(K_{2n})$ for $n \ge 3$, \emph{can} have their edge sets partitioned into copies of their cores.

%

\renewcommand{\bibname}{References}

\bibliographystyle{plain}


\end{document}